\documentclass{amsart}

\usepackage{amsfonts, amsmath, amssymb, amsthm, amsxtra, latexsym}

\theoremstyle{plain}

\theoremstyle{definition}

\newcommand{\R}{{\mathbb R}}
\newcommand{\C}{{\mathbb C}}

\numberwithin{equation}{section}

\begin{document}


\title[Remarks on another proof of Paley--Wiener theorems]{Some remarks on another proof of geometrical Paley--Wiener theorems for the Dunkl transform}


\translator{}

\dedicatory{}

\author{Marcel de Jeu}

\begin{abstract}
We argue that another proof by Trim\`eche of the geometrical form of the
Paley--Wiener theorems for the Dunkl transform is not correct.
\end{abstract}

\date{}

\subjclass[2000]{Primary 33C52; Secondary 33C67}

\keywords{Dunkl transform, Paley--Wiener theorem}

\address{M.F.E.~de~Jeu\\
         Mathematical Institute\\
         Leiden University\\
         P.O. Box 9512\\
         2300 RA Leiden\\
         The Netherlands}



\email{mdejeu@math.leidenuniv.nl}

\urladdr{}

\maketitle


\section{Introduction}
In \cite{trimeche} a proof of the geometrical form of the Paley--Wiener
theorems for the Dunkl transform is presented. It is argued in
\cite{pw_remarks} that this proof is not correct. A new proof has appeared in
\cite{trimeche_2}. In this note we argue that this new proof is not correct.

The material which is presented below has previously been communicated to the
author. It is our opinion that at this moment the geometrical forms of the
Paley--Wiener theorems for the Dunkl transform are still unproven, and that
only partial results have been established \cite{pw_new}.

\section{Arguments}

In \cite{trimeche_2} the geometric form of the Paley-Wiener theorem for the
Dunkl transform is stated for functions as Theorem~5.1. Our arguments concern
the proof of this theorem. In formulating them, we will use the notation and
definitions of \cite{trimeche_2}.

In the proof of Theorem~5.1, the function $\phi$ on $\R^d$ is defined in
equation (5.3) as
\begin{equation}\label{eq:phidefinition}
\phi(x)=\int_{\R^d} f(y)K(iy,x)\frac{\omega_k(y)}{(1+\Vert
y\Vert^2)^p}dy\quad(x\in\R^d).
\end{equation}
Here $f$ is a function of Paley--Wiener type corresponding to a $W$-invariant
compact convex subset $E$ of $\R^d$, as in the statement of Theorem~5.1, and
$p$ is an integer such that $p>\gamma+\frac{d}{2}+1$. Since the constant
$\gamma$ is assumed to be strictly positive in line~-10 of page 4, we see that
that $p\geq 2$.

After a computation involving Riemann sums and contour integration, it is then
concluded on line~9 of page~13 that $\phi$ has support in the set $E$. In
particular, $\phi$ has compact support. Since it has already been observed in
line~10 on page~11 that $\phi$ is smooth, the easy part of the Paley--Wiener
theorem implies that the Dunkl transform of $\phi$ is of Paley--Wiener type. In
particular, the Dunkl transform of $\phi$ has an entire extension to $\C^d$.
The inversion theorem and our equation~\eqref{eq:phidefinition} therefore show
that the function
\begin{equation}\label{eq:hasanextension}
\frac{f(y)}{(1+\Vert y\Vert^2)^p}
\end{equation}
on $\R^d$ has an entire extension to $\C^d$.

Now take $E$ to be an arbitrary closed ball, centered at the origin. Certainly,
if $f$ is the ordinary Fourier transform of a smooth function with support in
$E$, then $f$ is of Paley--Wiener type corresponding to the $W$-invariant
compact convex set $E$, so that the above reasoning applies to $f$.

Combining the previous two paragraphs, we conclude that the function in our
equation~\eqref{eq:hasanextension} has an entire extension from $\R^d$ to
$\C^d$, for all $f$ which are the ordinary Fourier transform of a compactly
supported smooth function. Since $f$ has then itself an entire extension, we
see that the meromorphic function
\[
\frac{f(y)}{(1+(y,y))^p}
\]
on $\C^d$ has a removable singularity along the divisor $\{y\in\C^d\mid
(y,y)=-1\}$, where $(\,.\,,\,.\,)$ is the standard bilinear form on $\C^d$.
Since $p\geq 2>0$, we conclude that each function $f$, which is the entire
extension of the ordinary Fourier transform of a smooth compactly supported
function, vanishes on this divisor.

Now let $f$ be the Fourier transform of a non-zero positive smooth compactly
supported function. Then $f(i,0,\ldots,0)>0$. Since
$((i,0,\ldots,),(i,0,\ldots,))=-1$ we must also have $f(i,0,\ldots,)=0$ by the
reasoning above. This is a contradiction.

\medskip
The above arguments are not based on the details of the computation with
Riemann sums and contour integration, but they are concerned with the
impossibility of the support of $\phi$ being contained in $E$. These arguments
therefore show not only that this computation contains a technical inaccuracy,
but they also show that the technique of this computation can not be corrected,
since the conclusion which the computation is supposed to yield does not hold.



\begin{thebibliography}{99}


\bibitem{pw_remarks}
M.F.E.~de~Jeu, \emph{Some remarks on a proof of geometrical Paley-Wiener
theorems for the Dunkl transform}, Preprint (2004). ArXiv: math.CA/0404293.

\bibitem{pw_new}
M.F.E.~de~Jeu, \emph{Paley--Wiener theorems for the Dunkl transform}, Preprint
(2004). ArXiv: math.CA/0404439.

\bibitem{trimeche}K.~Trim\`eche, \emph{Paley-Wiener theorems for the Dunkl transform and Dunkl translation operators}, Integral Transforms Spec.\ Funct.\  {\bf 13} (2002), 17--38.

\bibitem{trimeche_2} K.~Trim\`eche, \emph{Another proofs of the geometrical forms of Paley--Wiener theorems for the Dunkl transform and inversion formulas for the Dunkl interwining operator and for its
dual}, Preprint (2004). ArXiv: math.CA/0405050 v1.





\end{thebibliography}
\end{document}